\documentclass[11pt]{article}
\usepackage{amsmath}
\usepackage{amssymb}

\usepackage{theorem}
\numberwithin{equation}{section}

\newtheorem{theorem}{Theorem}[section]
\newtheorem{proposition}[theorem]{Proposition}
\newtheorem{corollary}[theorem]{Corollary}

\newtheorem{lemma}[theorem]{Lemma}
\newtheorem{definition}{Definition}[section]

\begin{document}
\title{Global well-posedness and scattering for the defocusing, mass - critical generalized KdV equation}
\date{\today}
\author{Benjamin Dodson}
\maketitle

\noindent \textbf{Abstract:} In this paper we prove that the defocusing, mass - critical generalized KdV initial value problem is globally well-posed and scattering for $u_{0} \in L^{2}(\mathbf{R})$. We prove this via a concentration compactness argument.

\section{Introduction}
In this paper we plan to study the global well - posedness theory for the initial value problem for the defocusing generalized KdV equation,

\begin{equation}\label{1.1}
 \partial_{t} u + \partial_{xxx} u = \partial_{x}(u^{5}), \hspace{5mm} u(0) \in L^{2}(\mathbf{R}), \hspace{5mm} x \in \mathbf{R}, t \in \mathbf{R}.
\end{equation}

\noindent The set of solutions of $(\ref{1.1})$ is invariant under the scaling

\begin{equation}\label{1.2}
 u_{\lambda}(x,t) = \lambda^{1/2} u(\lambda^{3} t, \lambda x)
\end{equation}

\noindent in the sense that if $u$ solves $(\ref{1.1})$ then so does $u_{\lambda}$ with initial datum

\begin{equation}\label{1.3}
 u_{\lambda}(0, x) = \lambda^{1/2} u(0,\lambda x).
\end{equation}

\noindent Notice that $\| u_{\lambda}(0,x) \|_{L^{2}(\mathbf{R})} = \| u(0,x) \|_{L^{2}(\mathbf{R})}$, so $(\ref{1.1})$ is an $L^{2}$ critical generalized KdV equation. The $L^{2}$ norm, or mass, is conserved under the flow $(\ref{1.1})$.

\begin{equation}\label{1.4}
 M(u(t)) = \int_{\mathbf{R}} |u(t,x)|^{2} dx = M(u(0)).
\end{equation}

\noindent Another conserved quantity of $(\ref{1.1})$ is the energy

\begin{equation}\label{1.5}
 E(u(t)) = \frac{1}{2} \int_{-\infty}^{\infty} u_{x}^{2}(t,x) dx + \frac{1}{6} \int_{-\infty}^{\infty} u^{6}(t,x) dx = E(u(0)).
\end{equation}

\noindent We define a solution of $(\ref{1.1})$ to be a strong solution.

\begin{definition}[Solution]\label{d1.1}
A function $u : I \times \mathbf{R} \rightarrow \mathbf{R}$ on a non - empty interval $0 \in I \subset \mathbf{R}$ is a (strong) solution to $(\ref{1.1})$ if it lies in the class $C_{t}^{0} L_{x}^{2}(J \times \mathbf{R}) \cap L_{x}^{5} L_{t}^{10}(J \times \mathbf{R})$ for any compact $J \subset I$, and obeys the Duhamel formula

\begin{equation}\label{1.6}
u(t) = e^{-t \partial_{x}^{3}} u_{0} + \int_{0}^{t} e^{-(t - \tau) \partial_{x}^{3}} \partial_{x}(u^{5}(\tau)) d\tau.
\end{equation}

\noindent We refer to the interval $I$ as the lifespan of $u$. We say that $u$ is a maximal lifespan solution if the solution cannot be extended to any strictly larger interval. We say that $u$ is a global solution if $I = \mathbf{R}$.
\end{definition}

\noindent \cite{KPV} developed a global in time theory for initial data small enough in $L_{x}^{2}(\mathbf{R})$. The results turn local for arbitrary data with the time of existence depending on the shape of the initial data $u_{0}$ not just its size. In particular, if $u_{0}$ is a little bit more regular than $L_{x}^{2}(\mathbf{R})$, say $u_{0} \in H_{x}^{s}(\mathbf{R})$ for some $s > 0$, then a solution to $(\ref{1.1})$ exists on a time interval $[0, T]$, $T(\| u_{0} \|_{H_{x}^{s}(\mathbf{R})}) > 0$. This implies that a solution to $(\ref{1.1})$ is global if $u_{0} \in H_{x}^{1}(\mathbf{R})$.\vspace{5mm}

\noindent From $(\ref{d1.1})$ we can see that it is important to analyze the scattering size.

\begin{definition}[Scattering size]\label{d1.2}
\begin{equation}\label{1.7}
 S_{I}(u) = \int_{\mathbf{R}} (\int_{I} |u(t,x)|^{10} dt)^{1/2} dx = \| u \|_{L_{x}^{5} L_{t}^{10}(I \times \mathbf{R})}^{5}.
\end{equation}
\end{definition}

\noindent Associated with the notion of a solution is a corresponding notion of blowup.

\begin{definition}[Blowup]\label{d1.3}
 We say that a solution $u$ to $(\ref{1.1})$ blows up forward in time if there exists $t_{1} \in I$ such that

\begin{equation}\label{1.8}
 S_{[t_{1}, \sup(I))}(u) = \infty.
\end{equation}

\noindent and that $u$ blows up backward in time if there exists a time $t_{1} \in I$ such that

\begin{equation}\label{1.9}
S_{(\inf(I), t_{1}]}(u) = \infty.
\end{equation}
\end{definition}

\noindent This precisely corresponds to the impossibility of continuing the solution (in the case of blowup in finite time) or failure to scatter (in the case of blowup in infinite time). We summarize the results of \cite{KPV} below.

\begin{theorem}[Local well - posedness]\label{t1.4}
Given $u_{0} \in L_{x}^{2}(\mathbf{R})$ and $t_{0} \in \mathbf{R}$, there exists a unique maximal lifespan solution $u$ to $(\ref{1.1})$ with $u(t_{0}) = u_{0}$. We will write $I$ for the maximal lifespan. This solution also has the following properties:\vspace{5mm}

1. (Local existence) I is an open neighborhood of $t_{0}$.\vspace{5mm}

2. (Blowup criterion) If $\sup(I)$ is finite then $u$ blows up forward in time. If $\inf(I)$ is finite then $u$ blows up backward in time.\vspace{5mm}

3. (Scattering) If $\sup(I) = +\infty$ and $u$ does not blow up forward in time, then $u$ scatters forward in time. That is, there exists a unique $u_{+} \in L_{x}^{2}(\mathbf{R})$ such that

\begin{equation}\label{1.10}
 \lim_{t \rightarrow +\infty} \| u(t) - e^{-t \partial_{x}^{3}} u_{+} \|_{L_{x}^{2}(\mathbf{R})} = 0.
\end{equation}

\noindent Conversely, given $u_{+} \in L_{x}^{2}(\mathbf{R})$ there is a unique solution to $(\ref{1.1})$ in a neighborhood of $\infty$ so that $(\ref{1.10})$ holds. One can define scattering backward in time in a completely analogous manner.\vspace{5mm}

4. (Small data global existence) If $M(u_{0})$ is sufficiently small then $u$ is a global solution which does not blow up either forward or backward in time. Indeed, in this case

\begin{equation}\label{1.11}
S_{\mathbf{R}}(u) \lesssim M(u)^{5/2}.
\end{equation}
\end{theorem}

\noindent \textbf{Remark:} See \cite{CaWe} for the analogous result for the nonlinear Schr{\"o}dinger equation. In this paper we will prove

\begin{theorem}[Spacetime bounds for the mass - critical gKdV]\label{t1.5}
The defocusing mass - critical gKdV problem $(\ref{1.1})$ is globally well - posed for arbitrary initial data $u_{0} \in L_{x}^{2}(\mathbf{R})$. Furthermore, the global solution satisfies the following spacetime bounds

\begin{equation}\label{1.12}
 \| u \|_{L_{x}^{5} L_{t}^{10}(\mathbf{R} \times \mathbf{R})} \leq C(M(u_{0})).
\end{equation}

\noindent The function $C : [0, \infty) \rightarrow [0, \infty)$.
\end{theorem}

\noindent \textbf{Remark:} This paper does not consider the focusing problem at all. See \cite{KPV1} and \cite{KKSV} for more information on this topic and the conjectured result.\vspace{5mm}

\noindent This theorem is proved using concentration compactness. \cite{KPV1} demonstrated that if a solution to $(\ref{1.1})$ blows up in finite time $T_{\ast} < \infty$, there exists a $C_{0}$ such that at least $C_{0}$ amount of mass must concentrate in a window of width $c(T_{\ast} - t)^{1/2} \| u(t) \|_{H_{x}^{s}}^{1/2s}$ for some $s > 0$.\vspace{5mm}

\noindent Later, \cite{KKSV} proved a conditional concentration compactness result.

\begin{theorem}[Concentration compactness theorem]\label{t1.6}
Assume that the defocusing mass - critical nonlinear Schr{\"o}dinger equation in one dimension,

\begin{equation}\label{1.13}
 (i \partial_{t} + \partial_{xx})v = |v|^{4} v
\end{equation}

\noindent has global spacetime bounds

\begin{equation}\label{1.14}
 \int_{\mathbf{R}} \int_{\mathbf{R}} |v(t,x)|^{6} dx dt \leq C(M(v(0,x))).
\end{equation}

\noindent Then if theorem $\ref{t1.5}$ fails to be true, there exists a critical mass $0 < M_{c} < \infty$ such that $u$ is a blowup solution in both time directions to $(\ref{1.1})$ on some maximal interval $I$, $M(u(t)) = M_{c}$, and $\{ u(t) : t \in I \} \subset \{ \lambda^{1/2} f(\lambda(x + x_{0})) : \lambda \in (0, \infty), x_{0} \in \mathbf{R}, f \in K \}$ for some compact $K \subset L_{x}^{2}(\mathbf{R})$.
\end{theorem}

\noindent Subsequently \cite{D3} proved that a solution to $(\ref{1.13})$ does have the global spacetime bounds $(\ref{1.14})$. Therefore, at this point it only remains to rule out the minimal mass blowup solution described in theorem $\ref{t1.6}$. Notice that modulo symmetries in $x_{0}$ and $\lambda$ the minimal mass blowup solution described in theorem $\ref{t1.6}$ lies in a precompact set. Therefore, a sequence of solutions will have a convergent subsequence modulo symmetries in $x_{0}$ and $\lambda$. For any $t \in I$ let $N(t) \in (0, \infty)$ and $x(t) \in \mathbf{R}$ be the scale function and spatial function respectively such that

\begin{equation}\label{1.15}
 N(t)^{-1/2} u(N(t)^{-1} (x - x(t))) \in K.
\end{equation}

\noindent \textbf{Remark:} We have some flexibility with regard to the $N(t)$, $x(t)$ and $K$ that we choose. This will be discussed in the concentration compactness section. To rule out the minimal mass blowup solution in theorem $\ref{t1.6}$ it suffices to rule out one of three scenarios,

\noindent 1. The self - similar scenario.\vspace{5mm}

\begin{equation}\label{1.16}
 N(t) \sim t^{-1/3}, \hspace{5mm} t \in (0, \infty)
\end{equation}

\noindent 2. The double rapid cascade.\vspace{5mm}

\begin{equation}\label{1.17}
 N(t) \geq 1, \hspace{5mm} N(0) = 1, \hspace{5mm} \int_{I} N(t)^{2} dt \lesssim 1,
\end{equation}

\begin{equation}\label{1.17.1}
 \lim_{t \nearrow \sup(I)} N(t) = \lim_{t \searrow \inf(I)} N(t) = +\infty.
\end{equation}

\noindent 3. The quasisoliton solution.\vspace{5mm}

\begin{equation}\label{1.18}
 \int_{J} N(t)^{3} dt \sim \mathcal J, \hspace{5mm} \int_{J} N(t)^{2} \lesssim \mathcal J,
\end{equation}

\begin{equation}\label{1.19}
 E(u(t)) \lesssim 1,
\end{equation}

\noindent for some $\mathcal J$ large, $J \subset I$.\vspace{5mm}

\noindent The first two scenarios are precluded by an additional regularity argument. We use concentration compactness to show that in cases one and two $E(u(t)) \lesssim 1$, which prevents $N(t) \nearrow \infty$.\vspace{5mm}

\noindent To rule out the quasisoliton we construct an interaction Morawetz estimate. We rely on the papers of \cite{TTao} and then \cite{KwSh}, which proved the nonexistence of a soliton solution to the generalized KdV equation by showing that the center of energy moves to the left faster than the center of mass. We utilize the computations in \cite{TTao} to produce an interaction Morawetz estimate that is similar in flavor to the interaction Morawetz estimate of \cite{D5}. This rules out the final scenario, proving theorem $\ref{t1.5}$.\vspace{5mm}

\noindent In section two we discuss some properties of the linear solution to the Airy equation $(\partial_{t} + \partial_{xxx})u = 0$ as well as estimates for the nonlinear equation $(\ref{1.1})$. Most of these estimates can be found in \cite{KPV} and \cite{KKSV}. We also will discuss the $U_{\partial_{x}^{3}}$ and $V_{\partial_{x}^{3}}$ spaces of \cite{HHK}.\vspace{5mm}

\noindent In section three we will discuss the local conservation of the quantities mass and energy. We will use the computations of \cite{TTao}.\vspace{5mm}

\noindent In section four we will describe the concentration compactness of \cite{KKSV}. We will then discuss our three minimal mass blowup scenarios.\vspace{5mm}

\noindent In section five we will rule out the self - similar blowup scenario.\vspace{5mm}

\noindent In section six we will rule out the double rapid cascade.\vspace{5mm}

\noindent In section seven we will rule out the quasi - soliton.\vspace{5mm}

\noindent \textbf{Acknowledgments:} At this time the author would like to thank Luis Vega for sending him a copy of \cite{KPV1} and encouraging to work on the KdV problem.

\section{Linear Estimates}
We are interested in the mixed norm spaces

\begin{equation}\label{2.1}
 L_{x}^{p} L_{t}^{q}(I \times \mathbf{R}) = \{ F(x,t) : (\int_{\mathbf{R}} (\int_{I} |F(x,t)|^{q} dt)^{p/q} dx)^{1/p} < +\infty \},
\end{equation}

\noindent and

\begin{equation}\label{2.2}
 L_{t}^{p} L_{x}^{q}(I \times \mathbf{R}) = \{ F(t,x) : (\int_{I} (\int_{\mathbf{R}} |F(t,x)|^{q} dt)^{p/q} dx)^{1/p} < +\infty \}.
\end{equation}

\begin{definition}\label{d2.1}
$(p, q, \alpha)$ is an admissible triple if

\begin{equation}\label{2.3}
 \frac{1}{p} + \frac{1}{2q} = \frac{1}{4}, \hspace{5mm} \alpha = \frac{2}{q} - \frac{1}{p}, \hspace{5mm} 1 \leq p, q \leq \infty, \hspace{5mm} -\frac{1}{4} \leq \alpha \leq 1.
\end{equation}

\noindent If $(p, q, \alpha)$ is an admissible triple denote $(p, q, \alpha) \in \mathcal A$.
\end{definition}

\begin{proposition}[Linear estimates]\label{p2.2}
Let $u$ be a solution of the initial value problem

\begin{equation}\label{2.4}
\aligned
 (\partial_{t} + \partial_{x}^{3})u &= F, \hspace{5mm} t \in I, x \in \mathbf{R}, \\
u(0, x) &= u_{0}.
\endaligned
\end{equation}

\noindent Then for any admissible triples $(p_{j}, q_{j}, \alpha_{j})$, $j = 1, 2$,

\begin{equation}\label{2.5}
 \| D_{x}^{\alpha_{1}} u \|_{L_{x}^{p_{1}} L_{t}^{q_{1}}(I \times \mathbf{R})} \lesssim \| u_{0} \|_{L^{2}(\mathbf{R})} + \| D_{x}^{-\alpha_{2}} F \|_{L_{x}^{p_{2}'} L_{t}^{q_{2}'}(I \times \mathbf{R})}.
\end{equation}
\end{proposition}

\noindent \emph{Proof:} This was proved in \cite{KPV1}. $\Box$\vspace{5mm}

\noindent Taking a cue from the analysis of the nonlinear Schr{\"o}dinger equation (see for example \cite{Tao}), consider the analogue of the Strichartz spaces in the gKdV case.

\begin{definition}\label{d2.3}
Let

\begin{equation}\label{2.6}
\| u \|_{S^{0}(I \times \mathbf{R})} = \sup_{(p,q,\alpha) \in \mathcal A} \| D_{x}^{\alpha} u \|_{L_{x}^{p} L_{t}^{q}(I \times \mathbf{R})}.
\end{equation}

\noindent Then let $N^{0}(I \times \mathbf{R})$ be the dual of $S^{0}(I \times \mathbf{R})$ with appropriate norm.

\begin{equation}\label{2.7}
 \| F \|_{N^{0}(I \times \mathbf{R})} = \inf_{F = F_{1} + F_{2}} \| D_{x}^{1/4} F_{1} \|_{L_{x}^{4/3} L_{t}^{1}(I \times \mathbf{R})} + \| D_{x}^{-1} F_{2} \|_{L_{x}^{1} L_{t}^{2}(I \times \mathbf{R})}.
\end{equation}
\end{definition}

\begin{lemma}[More linear estimates]\label{l2.4}
If $u$ is a solution to $(\ref{2.4})$ then

\begin{equation}\label{2.8}
 \| u \|_{S^{0}(I \times \mathbf{R})} + \| u \|_{L_{t}^{\infty} L_{x}^{2}(I \times \mathbf{R})} \lesssim \| u_{0} \|_{L_{x}^{2}(\mathbf{R})} + \| F_{1} \|_{N^{0}(I \times \mathbf{R})} + \| F_{2} \|_{L_{t}^{1} L_{x}^{2}(I \times \mathbf{R})},
\end{equation}

\noindent for any $F = F_{1} + F_{2}$ decomposition.

\end{lemma}

\noindent \emph{Proof:} See \cite{Kato}, \cite{KPV2}, \cite{KPV}, and \cite{KPV1}. $\Box$\vspace{5mm}

\noindent In this paper it is useful to use the $U_{\partial_{x}^{3}}^{2}$ and $V_{\partial_{x}^{3}}^{2}$ spaces of \cite{HHK}.

\begin{definition}\label{d2.5}
Let $1 \leq p < \infty$. $u$ is a $U_{\partial_{x}^{3}}^{p}(I \times \mathbf{R})$ atom if $[t_{0}, t_{1}]$, $[t_{1}, t_{2}]$, ... is a partition of $I$,

\begin{equation}\label{2.9}
 u = \sum_{t_{j} \nearrow} 1_{[t_{j}, t_{j + 1}]}(t) e^{-t \partial_{x}^{3}} u(t_{j}), 
\end{equation}

\begin{equation}\label{2.10}
 \sum_{t_{j} \nearrow} \| u(t_{j}) \|_{L_{x}^{2}(\mathbf{R})}^{p} = 1.
\end{equation}

\noindent Then define the norm

\begin{equation}\label{2.11}
 \| u \|_{U_{\partial_{x}^{3}}^{p}(I \times \mathbf{R})} = \inf \{ \sum_{\lambda} |c_{\lambda}| : \sum_{\lambda} c_{\lambda} u_{\lambda} = u \hspace{5mm} a.e., \hspace{5mm} u_{\lambda} \text{ is a } U_{\partial_{x}^{3}}^{p} \text{ atom } \}.
\end{equation}

\noindent Let

\begin{equation}\label{2.12}
 \| v \|_{V_{\partial_{x}^{3}}^{p}(I \times \mathbf{R})}^{p} = \sup_{\{ t_{j} \nearrow \}} \sum_{t_{j} \nearrow} \| e^{t_{j} \partial_{x}^{3}} v(t_{j}) - e^{t_{j + 1} \partial_{x}^{3}} v(t_{j + 1}) \|_{L_{x}^{2}(\mathbf{R})}^{p}.
\end{equation}

\begin{equation}\label{2.13}
 \| F \|_{DU_{\partial_{x}^{3}}^{p}(I \times \mathbf{R})} = \inf \{ \| u \|_{U_{\partial_{x}^{3}}^{p}(I \times \mathbf{R})} : (\partial_{t} + \partial_{x}^{3}) u = F \}.
\end{equation}
\end{definition}

\noindent \textbf{Remark:} By checking individual atoms and direct calculation, $U_{\partial_{x}^{3}}^{p} \subset U_{\partial_{x}^{3}}^{q}$, $V_{\partial_{x}^{3}}^{p} \subset V_{\partial_{x}^{3}}^{q}$ when $p < q$.\vspace{5mm}

\noindent \textbf{Remark:} By checking individual atoms,

\begin{equation}\label{2.13.1}
\| u \|_{S^{0}(I \times \mathbf{R})} \lesssim \| u \|_{U_{\partial_{x}^{3}}^{2}(I \times \mathbf{R})}.
\end{equation}

\noindent It can be verified by direct calculation (see \cite{HHK}) that

\begin{equation}\label{2.13.2}
\| F \|_{DU_{\partial_{x}^{3}}^{p}(I \times \mathbf{R})} = \sup_{\| v \|_{V_{\partial_{x}^{3}}^{p'}(I \times \mathbf{R})} = 1} \langle v, F \rangle.
\end{equation}

\begin{lemma}\label{l2.6}
For a decomposition $F = F_{1} + F_{2}$,

\begin{equation}\label{2.14}
\| F \|_{DU_{\partial_{x}^{3}}^{2}(I \times \mathbf{R})} \lesssim \| |\partial_{x}|^{-1/6} F_{1} \|_{L_{t,x}^{6/5}(I \times \mathbf{R})} + \| F_{2} \|_{L_{x}^{5/4} L_{t}^{10/9}(I \times \mathbf{R})}.
\end{equation}

\begin{equation}\label{2.15}
 \| \partial_{x}(u^{5}) \| _{DU_{\partial_{x}^{3}}^{2}(I \times \mathbf{R})} \lesssim \| u \|_{S^{0}(I \times \mathbf{R})}^{5}.
\end{equation}
\end{lemma}

\noindent \emph{Proof:} The first inequality follows from the embedding $V_{\partial_{x}^{3}}^{2} \subset U_{\partial_{x}^{3}}^{p}$ for any $p > 2$ (see \cite{HHK}). It can be verified by checking individual atoms that

\begin{equation}\label{2.16}
 \| |\partial_{x}|^{1/6} v \|_{L_{t,x}^{6}(I \times \mathbf{R})} + \| v \|_{L_{x}^{5} L_{t}^{10}(I \times \mathbf{R})} \lesssim \| v \|_{U_{\partial_{x}^{3}}^{5}(I \times \mathbf{R})} \lesssim \| v \|_{V_{\partial_{x}^{3}}^{2}(I \times \mathbf{R})} = 1.
\end{equation}

\noindent Next,

\begin{equation}\label{2.17}
\| \partial_{x}(u^{5}) \|_{L_{x}^{5/4} L_{t}^{10/9}(I \times \mathbf{R})} \lesssim \| \partial_{x} u \|_{L_{x}^{\infty} L_{t}^{2}(I \times \mathbf{R})} \| u \|_{L_{x}^{5} L_{t}^{10}(I \times \mathbf{R})}^{4} \leq \| u \|_{S^{0}(I \times \mathbf{R})}^{5}.
\end{equation}

\noindent This proves $(\ref{2.15})$. $\Box$\vspace{5mm}

\noindent We also make use of the dispersive estimate. 

\begin{lemma}[Dispersive estimate]\label{l2.6.1}
For $2 \leq p \leq \infty$,

\begin{equation}\label{2.17.1}
\| e^{-t \partial_{x}^{3}} u_{0} \|_{L_{x}^{p}(\mathbf{R})} \lesssim t^{-\frac{2}{3}(\frac{1}{2} - \frac{1}{p})} \| u_{0} \|_{L_{x}^{p'}(\mathbf{R})}.
\end{equation}
\end{lemma}

\noindent Finally it is useful to quote a long - time stability theorem.

\begin{theorem}[Long - time stability for the mass - critical gKdV]\label{t2.7}
Let $I$ be a time interval containing zero and let $\tilde{u}$ bea solution to 

\begin{equation}\label{2.18}
 (\partial_{t} + \partial_{xxx}) \tilde{u} = \partial_{x} (\tilde{u}^{5}) + e, \hspace{5mm} \tilde{u}(0,x) = \tilde{u}_{0}(x).
\end{equation}

\noindent Assume that

\begin{equation}\label{2.19}
 \| \tilde{u} \|_{L_{t}^{\infty} L_{x}^{2}(I \times \mathbf{R})} \leq M, \| \tilde{u} \|_{L_{x}^{5} L_{t}^{10}(I \times \mathbf{R})} \leq L
\end{equation}

\noindent for some positive constants $M$ and $L$. Let $u_{0}$ be such that

\begin{equation}\label{2.20}
 \| u_{0} - \tilde{u}_{0} \|_{L_{x}^{2}} \leq M'.
\end{equation}

\noindent Assume also the smallness conditions

\begin{equation}\label{2.21}
\aligned
\| e^{-t \partial_{x}^{3}} (u_{0} - \tilde{u}_{0}) \|_{L_{x}^{5} L_{t}^{10}(I \times \mathbf{R})} &\leq \epsilon, \\
\| e \|_{N^{0}(I \times \mathbf{R})} &\leq \epsilon,
\endaligned
\end{equation}

\noindent for some small $0 < \epsilon < \epsilon_{1}(M, M', L)$. Then there exists a solution $u$ to $(\ref{1.1})$ on $I \times \mathbf{R}$ with initial data $u_{0}$ at time $t = 0$ satisfying

\begin{equation}\label{2.22}
 \aligned
\| u - \tilde{u} \|_{L_{x}^{5} L_{t}^{10}(I \times \mathbf{R})} \leq C(M, M', L) \epsilon, \\
\| u^{5} - \tilde{u}^{5} \|_{L_{x}^{1} L_{t}^{2}(I \times \mathbf{R})} \leq C(M, M', L) \epsilon, \\
\| u - \tilde{u} \|_{L_{t}^{\infty} L_{x}^{2}(I \times \mathbf{R})} + \| u - \tilde{u} \|_{S^{0}(I \times \mathbf{R})} \leq C(M, M', L), \\
\| u \|_{L_{t}^{\infty} L_{x}^{2}(I \times \mathbf{R})} + \| u \|_{S^{0}(I \times \mathbf{R})} \leq C(M, M', L).
\endaligned
\end{equation}
\end{theorem}

\noindent \emph{Proof:} See \cite{KKSV}. $\Box$\vspace{5mm}

\noindent In particular, this theorem implies that if $u_{0}^{n} \rightarrow u_{0}$ strongly in $L^{2}$, and $u$ is the solution to $(\ref{1.1})$ on $I \subset \mathbf{R}$ with initial data $u_{0}$, then for any $J \subset I$,

\begin{equation}\label{2.23}
 \| u \|_{S^{0}(J \times \mathbf{R})} \leq C < \infty,
\end{equation}

\noindent then $u^{n} \rightarrow u$ in $S^{0}(I \times \mathbf{R})$ and $L_{t}^{\infty} L_{x}^{2}(I \times \mathbf{R})$, where $u^{n}$ is the solution to $(\ref{1.1})$ with initial data $u_{0}^{n}$.

\section{Local Conservation of mass and energy}
In this section we list the local conservation laws used in many places, for example \cite{TTao} and \cite{KwSh}.

\begin{definition}[Mass density and mass current]\label{d3.1}
The mass density is given by

\begin{equation}\label{3.1}
 \rho(t,x) = u(t,x)^{2}.
\end{equation}

\noindent The mass current is given by

\begin{equation}\label{3.2}
j(t,x) = 3 u_{x}(t,x)^{2} + \frac{5}{3} u(t,x)^{6}.
\end{equation}
\end{definition}

\begin{definition}[Energy density and energy current]\label{d3.2}
The energy density is given by

\begin{equation}\label{3.3}
e(t,x) = \frac{1}{2} u_{x}(t,x)^{2} + \frac{1}{6} u(t,x)^{6}.
\end{equation}

\noindent The energy current is given by

\begin{equation}\label{3.4}
k(t,x) = \frac{3}{2} u_{xx}(t,x)^{2} + 2 u(t,x))^{4} u_{x}(t,x)^{2} + \frac{1}{2} u(t,x)^{10}.
\end{equation}
\end{definition}

\noindent A routine computation verifies (for Schwartz solutions, at least) the pointwise conservation laws

\begin{equation}\label{3.5}
 \rho_{t} + \rho_{xxx} = j_{x},
\end{equation}

\begin{equation}\label{3.6}
 e_{t} + e_{xxx} = k_{x}.
\end{equation}

\noindent In section seven we will make use of the monotonicity formula.

\begin{lemma}[Monotonicity formula]\label{l3.3}
For a smooth function $u$,

\begin{equation}\label{3.7}
 (\int \rho(t,x) dx)(\int k(t,x) dx) - (\int e(t,x) dx)(\int j(t,x) dx) > 0.
\end{equation}

\end{lemma}

\noindent \emph{Proof:} See \cite{TTao}. $\Box$\vspace{5mm}

\noindent \textbf{Remark:} Frequently in this paper it will be necessary to integrate by parts. This paper will always assume that the solution is smooth in space and time. An arbitrary solution can be well approximated by a smooth solution, and the bounds obtained will not depend on the smoothness of $u$. Similar computations are done in the case of the interaction Morawetz estimate for the Schr{\"o}dinger equation. See for example \cite{CKSTT2}.

\section{Concentration Compactness}

\noindent An important step in the study of the mass critical generalized KdV was the reduction of \cite{KKSV} to solutions that are almost periodic modulo symmetries.

\begin{definition}[Almost periodic modulo symmetries]\label{d4.1}
A solution $u$ to (the mKdV problem) with lifespan $I$ is said to be almost periodic modulo symmetries if there exist functions $N : I \rightarrow \mathbf{R}^{+}$, $x : I \rightarrow \mathbf{R}$, $C : \mathbf{R}^{+} \rightarrow \mathbf{R}^{+}$ such that for all $t \in I$ and $\eta > 0$,

\begin{equation}\label{4.1}
 \int_{|x - x(t)| \geq \frac{C(\eta)}{N(t)}} |u(t,x)|^{2} dx + \int_{|\xi| \geq C(\eta) N(t)} |\hat{u}(t,\xi)|^{2} d\xi < \eta.
\end{equation}

\noindent $N$ will be called the frequency scale function for a solution $u$, $x$ the spatial center function, and $C$ the compactness modulus function.
\end{definition}

\noindent \textbf{Remark:} The parameter $N(t)$ measures the frequency scale of the solution at time $t$, while $\frac{1}{N(t)}$ measures the spatial scale. We can multiply $N(t)$ by any function $\alpha(t)$, $0 < \epsilon < \alpha(t) < \frac{1}{\epsilon}$, provided we also modify the compactness modulus function accordingly.

\begin{theorem}[Arzela - Ascoli theorem]\label{t4.2}
A family of functions is precompact in $L_{x}^{2}(\mathbf{R})$ if and only if it is norm bounded and there exists a compactness modulus function $C$ such that

\begin{equation}\label{4.2}
\int_{|x| \geq C(\eta)} |f(x)|^{2} dx + \int_{|\xi| \geq C(\eta)} |\hat{f}(\xi)|^{2} d\xi \leq \eta
\end{equation}

\noindent for all functions $f$ in the family.
\end{theorem}

\noindent This implies that $f$ is almost periodic modulo symmetries if and only if for some compact subset $K \subset L_{x}^{2}(\mathbf{R})$,

\begin{equation}\label{4.3}
 \{ u(t) : t \in I \} \subseteq \{ \lambda^{1/2} f(\lambda (x + x_{0})) : \lambda \in (0, \infty), x_{0} \in \mathbf{R}, f \in K \}.
\end{equation}

\noindent Let

\begin{equation}\label{4.4}
 L(M) = \sup \{ S_{I}(u) : u : I \times \mathbf{R} \rightarrow \mathbf{R}, M(u) \leq M_{c} \}.
\end{equation}

\noindent The supremum is taken over all solutions $u : I \times \mathbf{R} \rightarrow \mathbf{R}$ obeying $M(u) \leq M$. For $M$ small, a small data result implies $L(M) \lesssim M^{5/2}$. This fact combined with theorem $\ref{t2.7}$ implies that failure of theorem $\ref{t1.5}$ is equivalent to the existence of a critical mass $M_{c} \in (0, \infty)$ such that

\begin{equation}\label{4.5}
 L(M) < \infty \text{ for } M < M_{c}, \hspace{5mm} L(M) = \infty \text{ for } M \geq M_{c},
\end{equation}

\begin{theorem}\label{t4.3}
 Assume theorem $\ref{t2.7}$ fails. Let $M_{c}$ denote the critical mass. Then there exists a maximal lifespan solution to the mass - critical gKdV with mass $M(u) = M_{c}$ which is almost periodic modulo symmetries and blows up both forward and backward in time. Also, $[0, \infty) \subset I$, $N(t) \leq 1$ for $t \geq 0$, and

\begin{equation}\label{4.6}
 |N'(t)| \lesssim N(t)^{4}, \hspace{5mm} |x'(t)| \lesssim N(t)^{2}.
\end{equation}

\noindent Moreover, there exists $\delta(u) > 0$ such that for any $t_{0} \in I$,

\begin{equation}\label{4.7}
 \| u \|_{S^{0}([t_{0}, t_{0} + \frac{\delta}{N(t_{0})^{3}}] \times \mathbf{R})} \lesssim 1.
\end{equation}
\end{theorem}

\noindent \emph{Proof:} See \cite{KKSV}. The proof of theorem $\ref{t4.3}$ was conditional on the assumption that the following mass - critical nonlinear Schr{\"o}dinger equation result was true. $\Box$\vspace{5mm}

\begin{lemma}[No waste lemma]\label{l4.3.1}
If $u$ is a minimal mass blowup solution to $(\ref{1.1})$ then for any $t \in I$,

\begin{equation}\label{4.7.1}
u(t) = \lim_{T \rightarrow \sup(I)} \int_{t}^{T} e^{-(t - \tau) \partial_{x}^{3}} \partial_{x}(u^{5})(\tau) d\tau = \lim_{T \rightarrow \inf(I)} \int_{t}^{T} e^{-(t - \tau) \partial_{x}^{3}} \partial_{x}(u^{5})(\tau) d\tau,
\end{equation}

\noindent weakly in $L_{x}^{2}(\mathbf{R})$.
\end{lemma}

\noindent \emph{Proof:} This follows in a similar manner to \cite{TVZ1}. If $\sup(T) = +\infty$ then $N(t) \rightarrow +\infty$ combined with $(\ref{4.1})$ implies

\begin{equation}\label{4.7.2}
\lim_{T \rightarrow \sup(I)} \langle e^{-(t - T) \partial_{x}^{3}} u(T), u(t) \rangle = 0.
\end{equation}

\noindent The same would be true if $N(T) \rightarrow 0$. If $N(T) \sim N(t)$ as $T \rightarrow \sup(I)$ then $\sup(I) = +\infty$. The dispersive estimate $(\ref{2.17.1})$ combined with $(\ref{4.1})$ implies that in this case also

\begin{equation}\label{4.7.3}
\lim_{T \rightarrow \sup(I)} \langle e^{-(t - T) \partial_{x}^{3}} u(T), u(t) \rangle = 0.
\end{equation}

\noindent $\Box$

\begin{theorem}\label{t4.4}
If $u$ is a solution to the one dimensional, mass - critical nonlinear Schr{\"o}dinger equation

\begin{equation}\label{4.8}
 (i \partial_{t} + \partial_{xx})u = |u|^{4} u,
\end{equation}

\noindent Then

\begin{equation}\label{4.9}
 \| u \|_{L_{t,x}^{6}(\mathbf{R} \times \mathbf{R})} \leq C(\| u(0,\cdot) \|_{L^{2}}).
\end{equation}
\end{theorem}

\noindent \emph{Proof:} See \cite{D3}. $\Box$\vspace{5mm}

\noindent \textbf{Remark:} At this point we will select one minimal mass blowup solution in the form of theorem $\ref{t4.3}$ and then show that this solution does not exist. Therefore we can abbreviate $A \leq C(u) B$ as $A \lesssim B$.\vspace{5mm}

\noindent We rule out three separate scenarios. Let

\begin{equation}\label{4.10}
t_{0}(T) = \inf \{ t \in [0, T] : N(t) = \inf_{t \in [0, T]} N(t) \}
\end{equation}

\noindent $N(t)$ attains its infimum on $[0, T]$ since $N(t)$ is continuous.\vspace{5mm}

\noindent \textbf{Case 1:} Self - similar solution.

\begin{equation}\label{4.11}
 \limsup_{T \rightarrow \infty} (\inf_{t \in [0, T]} N(t)) \cdot (\int_{0}^{T} N(t)^{2} dt) \leq C < +\infty.
\end{equation}

\begin{equation}\label{4.12}
 \limsup_{T \rightarrow \infty} \frac{\sup_{t \in [t_{0}, T]} N(t)}{N(t_{0}(T))} \leq C < +\infty.
\end{equation}

\noindent \textbf{Case 2:} Rapid double cascade.

\begin{equation}\label{4.13}
 \limsup_{T \rightarrow \infty} (\inf_{t \in [0, T]} N(t)) \cdot (\int_{0}^{T} N(t)^{2} dt) = C < \infty.
\end{equation}

\begin{equation}\label{4.14}
 \limsup_{T \rightarrow \infty} \frac{\sup_{t \in [t_{0}(T), T]} N(t)}{N(t_{0}(T))} = +\infty.
\end{equation}

\noindent \textbf{Case 3:} Quasi - soliton.

\begin{equation}\label{4.15}
 \limsup_{T \rightarrow \infty} (\inf_{t \in [0, T]} N(t)) (\int_{0}^{T} N(t)^{2} dt) = +\infty.
\end{equation}

\section{Self - Similar solution}

\noindent $(\ref{4.11})$ implies

\begin{equation}\label{5.1}
\liminf_{t \rightarrow \infty} N(t) = 0.
\end{equation}

\noindent Then $(\ref{4.12})$ implies that $N(t) \rightarrow 0$ as $t \rightarrow \infty$. For any integer $l \geq 0$ let

\begin{equation}\label{5.2}
 t_{l} = \inf \{ t : N(t) = 2^{-l} \}.
\end{equation}

\noindent Clearly $t_{0} = 1$. By $(\ref{4.11})$

\begin{equation}\label{5.3}
 2^{-l} t_{l} 2^{-2l} \leq C,
\end{equation}

\noindent so for any $l$, $t_{l} \lesssim 2^{3l}$. On the other hand $|N'(t)| \lesssim N(t)^{4}$ and $(\ref{4.12})$ imply

\begin{equation}\label{5.4}
 2^{-l} \leq \int_{t_{l - 1}}^{t_{l}} |N'(t)| dt \lesssim 2^{-4l} (t_{l} - t_{l - 1}) \leq 2^{-4l} t_{l}.
\end{equation}

\noindent This implies $t_{l} \gtrsim 2^{3l}$ and therefore $t_{l} \sim 2^{3l}$, so for $t \geq 1$, $(\ref{4.12})$ implies that $N(t) \sim t^{-1/3}$. Possibly after modifying $C(\eta)$ by a constant, let $N(t) = 1$ for $t \in [0, 1]$, $N(t) = t^{-1/3}$ for $t \in [1, \infty)$.\vspace{5mm}

\noindent Let $x(0) = 0$. $|x'(t) \lesssim N(t)^{2}$ so $|x(t)| \lesssim t^{1/3}$. Therefore, again after modifying $C(\eta)$ by a constant, for any $\eta > 0$ there exists $C(\eta) < \infty$ such that

\begin{equation}\label{5.5}
 \int_{|x| \geq \frac{C(\eta)}{N(t)}} u(t,x)^{2} dx + \int_{|\xi| \geq C(\eta) N(t)} |\hat{u}(t,\xi)|^{2} d\xi < \eta.
\end{equation}

\noindent Now take a sequence $t_{n} \rightarrow +\infty$ and let

\begin{equation}\label{5.6}
 u_{0}^{n} = \frac{1}{N(t_{n})^{1/3}} u(\frac{x}{N(t_{n})}).
\end{equation}

\noindent Then, passing to a subsequence, $u_{0}^{n} \rightarrow u_{0}$ in $L^{2}$ and if $u(1,\cdot) = u_{0}(\cdot)$, $u$ solves the mass critical mKdV, then $u$ is a self - similar blowup solution on $(0, \infty)$ and $N(t) = t^{-1/3}$. We then prove

\begin{theorem}[Additional regularity]\label{t5.1}
 If $u$ is a self - similar solution to the mass critical gKdV equation then $u(1) \in H_{x}^{1}(\mathbf{R}) \cap L^{6}(\mathbf{R})$.
\end{theorem}

\begin{corollary}[No self - similar solution]\label{c5.2}
 There does not exist a self - similar solution.
\end{corollary}

\emph{Proof:} conservation of energy contradicts $N(t) \rightarrow +\infty$ as $t \rightarrow 0$.

\noindent \emph{Proof of theorem $\ref{t5.1}$:} This proof is very similar to the additional regularity proof in \cite{KVZ}, \cite{KTV}, and \cite{TVZ2} for the self - similar blowup solution for the nonlinear Schr{\"o}dinger equation. The proof has two steps. First, using the double Duhamel formula we prove that a self - similar solution must possess some additional regularity. More precisely, for some $s > 0$,

\begin{equation}\label{5.6.1}
\| u \|_{H_{x}^{s}(\mathbf{R})} \sim t^{-s/3}.
\end{equation}

\noindent Then we argue by induction to show that in fact $u \in H_{x}^{1}(\mathbf{R})$. Let

\begin{equation}\label{5.7}
 \mathcal M(A) = \sup_{T} \| u_{\geq A T^{-1/3}} \|_{L_{t}^{\infty} L_{x}^{2}([T, 2T] \times \mathbf{R})},
\end{equation}

\begin{equation}\label{5.8}
 \mathcal S(A) = \sup_{T} \| u_{\geq A T^{-1/3}} \|_{U_{\partial_{x}^{3}}([T, 2T] \times \mathbf{R})},
\end{equation}

\begin{equation}\label{5.9}
 \mathcal N(A) = \sup_{T} \| P_{\geq A T^{-1/3}} \partial_{x}(u^{5}) \|_{U_{\partial_{x}^{3}}([T, 2T] \times \mathbf{R})}.
\end{equation}

\noindent By Duhamel's principle,

\begin{equation}\label{5.10}
 \mathcal S(A) \lesssim \mathcal M(A) + \mathcal N(A).
\end{equation}

\noindent Compactness in $L^{2}$ norm combined with the above estimate implies

\begin{equation}\label{5.11}
\lim_{A \rightarrow \infty} \mathcal M(A) + \mathcal S(A) + \mathcal N(A) = 0.
\end{equation}

\noindent Let $\alpha(k)$ be a frequency envelope that bounds $\| P_{2^{k}} u(1) \|_{L^{2}}$. Set $\delta = \frac{1}{40}$. Let

\begin{equation}\label{5.12}
 \alpha(k) = \sum_{j} 2^{-\delta |j - k|} \| P_{2^{j}} u(1) \|_{L^{2}}.
\end{equation}

\noindent Choose $\epsilon > 0$ very small, $k_{0}(\epsilon)$ sufficiently large so that

\begin{equation}\label{5.13}
 \mathcal M(2^{k_{0}/2}) + \mathcal S(2^{k_{0}/2}) + \mathcal N(2^{k_{0}/2}) < \epsilon,
\end{equation}

\begin{equation}\label{5.14}
 2^{-k_{0}} < \epsilon^{200},
\end{equation}

\begin{equation}\label{5.15}
 \sum_{k > k_{0}/2} \alpha(k)^{2} \leq \epsilon^{2}.
\end{equation}

\begin{theorem}\label{t5.3}
 For $k \geq k_{0}$,

\begin{equation}\label{5.16}
 \| P_{2^{k}} u \|_{U_{\partial_{x}^{3}}([1, 2^{6(k - k_{0})}] \times \mathbf{R})} \lesssim \alpha(k),
\end{equation}

\noindent and for $j > 6(k - k_{0})$,

\begin{equation}\label{5.17}
 \| P_{2^{k}} u \|_{U_{\partial_{x}^{3}}([2^{j}, 2^{j + 1}] \times \mathbf{R})} \lesssim 2^{\frac{1}{10}(j - 6(k - k_{0}))} \alpha(k).
\end{equation}
\end{theorem}

\noindent \emph{Proof:} We prove this by Duhamel's principle.

\begin{equation}\label{5.18}
 u(t) = e^{-(t - 1) \partial_{x}^{3}} u(1) + \int_{1}^{t} e^{-(t - \tau) \partial_{x}^{3}} \partial_{x}(u^{5}) d\tau.
\end{equation}

\begin{equation}\label{5.19}
 \| e^{-(t - 1) \partial_{x}^{3}} P_{2^{k}} u(1) \|_{U_{\partial_{x}^{3}}([1, 2^{6(k - k_{0})}] \times \mathbf{R})} \lesssim \alpha(k).
\end{equation}

\begin{equation}\label{5.20}
 \| P_{2^{k}} \partial_{x}(u^{5}) \|_{DU_{\partial_{x}^{3}}([1, 2^{6(k - k_{0})}] \times \mathbf{R})} \lesssim 2^{5k/6} \sum_{k_{1} \geq k} \| P_{k_{1}} u \|_{L_{t,x}^{6}([1, 2^{6(k - k_{0})}] \times \mathbf{R})}^{5}
\end{equation}

\begin{equation}\label{5.21}
 + 2^{5k/6} \| P_{k - 5 \leq \cdot \leq k + 5} u \|_{L_{x}^{\infty} L_{t}^{2}([1, 2^{6(k - k_{0} - 5)}] \times \mathbf{R})} \| P_{\leq k} u \|_{L_{x}^{24/5} L_{t}^{12}([1, 2^{6(k - k_{0} - 5)}] \times \mathbf{R})}^{4}
\end{equation}

\begin{equation}\label{5.22}
 + 2^{5k/6} \| P_{k - 5 \leq \cdot \leq k + 5} u \|_{L_{x}^{\infty} L_{t}^{2}([2^{6(k - k_{0} - 5)}, 2^{6(k - k_{0})}] \times \mathbf{R})} \| P_{\leq k} u \|_{L_{x}^{24/5} L_{t}^{12}([2^{6(k - k_{0} - 5)}, 2^{6(k - k_{0})}] \times \mathbf{R})}^{4}.
\end{equation}

\noindent By the local smoothing estimates and the concentration compactness result, for $j \geq k_{0}$,

\begin{equation}\label{5.23}
 \| P_{2^{j}} u \|_{L_{x}^{24/5} L_{t}^{12}([1, 2^{6(k - k_{0})}] \times \mathbf{R}} \lesssim 2^{j/24} \| P_{2^{j}} u \|_{U_{\partial_{x}^{3}}([1, 2^{6(j - k_{0})}] \times \mathbf{R})} + 2^{j/24} \| P_{2^{j}} u \|_{U_{\partial_{x}^{3}}([2^{6(j - k_{0})}, 2^{6(k - k_{0})}] \times \mathbf{R})}
\end{equation}

\begin{equation}\label{5.24}
 \lesssim \alpha(j) 2^{j/24} + 2^{j/24} (k - j) \epsilon.
\end{equation}

\noindent For $k_{0}/2 \leq j \leq k_{0}$,

\begin{equation}\label{5.25}
\| P_{2^{j}} u \|_{L_{x}^{24/5} L_{t}^{12}([1, 2^{6(k - k_{0})}] \times \mathbf{R}} \lesssim 2^{j/24} \epsilon(k - k_{0})^{5/24}.
\end{equation}

\noindent Finally for $j \leq k_{0}/2$,

\begin{equation}\label{5.26}
\| P_{2^{j}} u \|_{L_{x}^{24/5} L_{t}^{12}([1, 2^{6(k - k_{0})}] \times \mathbf{R}} \lesssim 2^{j/24} (k - k_{0})^{5/24}.
\end{equation}

\noindent Putting this all together,

\begin{equation}\label{5.27}
 2^{5k/6} \| P_{k - 5 \leq \cdot \leq k + 5} u \|_{L_{x}^{\infty} L_{t}^{2}([1, 2^{6(k - k_{0} - 5)}] \times \mathbf{R})} \| P_{\leq k} u \|_{L_{x}^{24/5} L_{t}^{12}([1, 2^{6(k - k_{0} - 5)}] \times \mathbf{R})}^{4} 
\end{equation}

\begin{equation}\label{5.28}
\aligned
\lesssim  2^{-k/6} \alpha(k) (\sum_{j \leq k} 2^{j/24} (\alpha(j) + \epsilon(k - j))^{4} + 2^{-k/6} \alpha(k) (\sum_{j \leq k_{0}} 2^{j/24} \epsilon (k - k_{0})^{5/24})^{4} \\ + 2^{-k/6} \alpha(k) (\sum_{j \leq k_{0}/2} 2^{j/24} (k - k_{0})^{5/24})^{4} \lesssim \alpha(k) \epsilon^{4}.
\endaligned
\end{equation}

\noindent Similarly,

\begin{equation}\label{5.29}
2^{5k/6} \| P_{k - 5 \leq \cdot \leq k + 5} u \|_{L_{x}^{\infty} L_{t}^{2}([2^{6(k - k_{0} - 5)}, 2^{6(k - k_{0})}] \times \mathbf{R})} \| P_{\leq k} u \|_{L_{x}^{24/5} L_{t}^{12}([2^{6(k - k_{0} - 5)}, 2^{6(k - k_{0})}] \times \mathbf{R})}^{4} \lesssim  \alpha(k) \epsilon^{4}.
\end{equation}

\noindent Finally,

\begin{equation}\label{5.30}
2^{5k/6} \sum_{k_{1} \geq k} \| P_{k_{1}} u \|_{L_{t,x}^{6}([1, 2^{6(k - k_{0})}] \times \mathbf{R})}^{5} \lesssim 2^{5k/6} \sum_{k_{1} \geq k} \alpha(k_{1})^{5} 2^{-5k_{1}/6} \lesssim \alpha(k) \epsilon^{4}.
\end{equation}

\noindent Now take $j > 6(k - k_{0})$.

\begin{equation}\label{5.31}
 \| P_{k} u \|_{U_{\partial_{x}^{3}}([2^{j}, 2^{j + 1}] \times \mathbf{R})} \leq \| P_{k} u(2^{j}) \|_{L^{2}} + C 2^{k} \| P_{k}(u^{5}) \|_{DU_{\partial_{x}^{3}}([2^{j}, 2^{j + 1}] \times \mathbf{R})}
\end{equation}

\begin{equation}\label{5.32}
  \leq \| P_{k} u \|_{U_{\partial_{x}^{3}}([2^{j - 1}, 2^{j}] \times \mathbf{R})} + C 2^{k} \| P_{k} (u^{5}) \|_{DU_{\partial_{x}^{3}}([2^{j}, 2^{j + 1}] \times \mathbf{R})}.
\end{equation}

\noindent By the same analysis as before,

\begin{equation}\label{5.33}
 2^{k} \| P_{k - 5 \leq \cdot \leq k + 5} u \|_{L_{x}^{\infty} L_{t}^{2}} \| P_{\leq k} u \|_{L_{x}^{24/5} L_{t}^{12}}^{4} + 2^{5k/6} \sum_{k_{1} > k} \| P_{k_{1}} u \|_{L_{t,x}^{6}}^{5}
\end{equation}

\begin{equation}\label{5.34}
 \lesssim 2^{(j - 6(k - k_{0}))/10} \alpha(k) \epsilon^{4} + 2^{5k/6} \epsilon^{4} \sum_{k \leq k_{1}} 2^{(j - 6(k - k_{0})/10} 2^{-5 k_{1}/6} \alpha(k_{1}) \lesssim 2^{(j - 6(k - k_{0}))/10} \alpha(2^{k}) \epsilon^{4}.
\end{equation}

\noindent Now make a bootstrapping argument. Let $A$ be the set of $T \in [1, \infty]$ such that for a large, fixed constant $C$,

\begin{equation}\label{5.35}
 \| P_{2^{k}} u \|_{U_{\partial_{x}^{3}}([1, 2^{6(k - k_{0})}] \cap [1, T] \times \mathbf{R})} \leq \frac{C}{2} \alpha(k),
\end{equation}

\noindent and for $j > 6(k - k_{0})$,

\begin{equation}\label{5.36}
 \| P_{2^{k}} u \|_{U_{\partial_{x}^{3}}([2^{j}, 2^{j + 1}] \cap [1, T] \times \mathbf{R})} \leq \frac{C}{2} 2^{\frac{1}{10}(j - 6(k - k_{0}))} \alpha(k).
\end{equation}

\noindent The set $A$ is nonempty since $1 \in A$, and is closed. It remains to show that $A$ is open. Suppose $A = [1, T_{0}]$. Then there exists $T_{0} < T < 2T_{0}$ such that

\begin{equation}\label{5.37}
 \| P_{2^{k}} u \|_{U_{\partial_{x}^{3}}([1, 2^{6(k - k_{0})}] \cap [1, T] \times \mathbf{R})} \leq C \alpha(k),
\end{equation}

\noindent and for $j > 6(k - k_{0})$,

\begin{equation}\label{5.38}
 \| P_{2^{k}} u \|_{U_{\partial_{x}^{3}}([2^{j}, 2^{j + 1}] \cap [1, T] \times \mathbf{R})} \leq C 2^{\frac{1}{10}(j - 6(k - k_{0}))} \alpha(k).
\end{equation}

\noindent For $\epsilon > 0$ sufficiently small,

\begin{equation}\label{5.39}
 \| P_{2^{k}} u \|_{U_{\partial_{x}^{3}}([1, 2^{6(k - k_{0})}] \cap [1, T] \times \mathbf{R})} \lesssim \alpha(k) + \alpha(k) \epsilon^{4},
\end{equation}

\noindent and for $j > 6(k - k_{0})$,

\begin{equation}\label{5.40}
 \| P_{2^{k}} u \|_{U_{\partial_{x}^{3}}([2^{j}, 2^{j + 1}] \cap [1, T] \times \mathbf{R})} \lesssim C 2^{\frac{1}{10}(j - 6(k - k_{0}))} \alpha(k) \epsilon^{4} + \frac{C}{2} 2^{\frac{1}{10}(j - 1 - 6(k - k_{0}))} \alpha(k).
\end{equation}

\noindent Choosing $\epsilon > 0$ sufficiently small, $C$ sufficiently large implies that the bounds for $C$ imply the bounds for $\frac{C}{2}$, which closes the bootstrap, proving that $A = [1, \infty)$. $\Box$\vspace{5mm}

\noindent Theorem $\ref{t5.3}$ implies that for $k > 6 k_{0}$,

\begin{equation}\label{5.43}
 \| P_{> 2^{k}} u \|_{U_{\partial_{x}^{3}}([2^{-5k/2}, 1] \times \mathbf{R})}^{2} \lesssim \sum_{j \geq k} (\sum_{j_{1}} 2^{-\delta |j_{1} - j|} \| P_{j_{1}} u(2^{-5k/2}) \|_{L_{x}^{2}(\mathbf{R})})^{2}
\end{equation}

\begin{equation}\label{5.44}
\aligned
 \lesssim \sum_{j_{1} \geq k - \frac{k_{0}}{2}} \| P_{j_{1}} u(2^{-5k/2}) \|_{L_{x}^{2}(\mathbf{R})}^{2} \sum_{j \geq k} 2^{-\delta|j_{1} - j|} + \sum_{j_{1} \leq k - \frac{k_{0}}{2}} \| P_{j_{1}} u(2^{-5k/2}) \|_{L_{x}^{2}(\mathbf{R})}^{2} \sum_{j \geq k} 2^{-\delta|j_{1} - j|} \\ \lesssim \epsilon^{2} + \mathcal M(2^{k_{0}/2})^{2} \lesssim \epsilon^{2}.
\endaligned
\end{equation}

\noindent By conservation of mass and the conditions on $k_{0}$.\vspace{5mm}

\noindent Another useful fact about self - similar solutions is that a self - similar solution rescales to another self - similar solution. The scaling

\begin{equation}\label{5.41}
 u(t,x) \mapsto \frac{1}{\lambda^{1/2}} u(\frac{t}{\lambda^{3}}, \frac{x}{\lambda}) = u_{\lambda}(t,x)
\end{equation}

\noindent with $\lambda = 2^{k}$ rescales the self - similar solution to a new self - similar solution with

\begin{equation}\label{5.42}
 u_{\lambda}(1) = \frac{1}{\lambda^{1/2}} u(\frac{1}{\lambda^{3}}, \frac{x}{\lambda}).
\end{equation}

\noindent The no - waste Duhamel formula $(\ref{4.7.1})$ gives the double Duhamel formula

\begin{equation}\label{5.45}
 \| P_{2^{k}} u(1) \|_{L^{2}}^{2} = \int_{0}^{1} \int_{1}^{\infty} \langle e^{(t - \tau) \partial_{x}^{3}} P_{k} \partial_{x}(u^{5}), P_{k} \partial_{x}(u^{5}) \rangle dt d\tau.
\end{equation}

\begin{equation}\label{5.46}
 \int_{2^{-5k/2}}^{1} \int_{1}^{2^{6(k - k_{0})}} \langle e^{(t - \tau) \partial_{x}^{3}} P_{k} \partial_{x}(u^{5}), P_{k} \partial_{x}(u^{5}) \rangle d\tau
\end{equation}

\begin{equation}\label{5.47}
 \lesssim 2^{5k/3} \| P_{k}(u^{5}) \|_{L_{t,x}^{6/5}([1, 2^{6(k - k_{0})}] \times \mathbf{R})} \| P_{k}(u^{5}) \|_{L_{t,x}^{6/5}([2^{-5k/2}, 1] \times \mathbf{R})}.
\end{equation}

\begin{equation}\label{5.48}
 2^{5k/6} \| P_{k}(u^{5}) \|_{L_{t,x}^{6/5}([1, 2^{6(k - k_{0})}] \times \mathbf{R})} \lesssim 2^{5k/6}(\sum_{j \geq k} \| P_{j} u \|_{L_{t,x}^{6}([1, 2^{6(k - k_{0})}] \times \mathbf{R})}^{5})
\end{equation}

\begin{equation}\label{5.49}
 + 2^{-k/6} \| P_{k - 5 \leq \cdot \leq k + 5} u \|_{L_{x}^{\infty} L_{t}^{2}([1, 2^{6(k - k_{0})}] \times \mathbf{R})} (\sum_{j \leq k} \| P_{j} u \|_{L_{x}^{24/5} L_{t}^{12}([1, 2^{6(k - k_{0})}] \times \mathbf{R})})^{4}.
\end{equation}

\noindent For all $j \geq k_{0}$, theorem $\ref{t5.3}$ implies

\begin{equation}\label{5.50}
 \| P_{j} u \|_{L_{x}^{24/5} L_{t}^{12}([1, 2^{6(k - k_{0})}] \times \mathbf{R})} \lesssim \alpha(j) 2^{\frac{(k - j)}{10}} 2^{j/24},
\end{equation}

\begin{equation}\label{5.51}
 \| P_{j} u \|_{L_{x}^{24/5} L_{t}^{12}([1, 2^{6(k - k_{0})}] \times \mathbf{R})} \lesssim \alpha(j) 2^{j/24} + 2^{j/24} (k - j)^{5/24},
\end{equation}

\noindent and for all $j$,

\begin{equation}\label{5.52}
 \| P_{j} u \|_{L_{x}^{24/5} L_{t}^{12}([1, 2^{6(k - k_{0})}] \times \mathbf{R})} \lesssim 2^{j/24} (k - k_{0})^{5/24}.
\end{equation}

\noindent Therefore,

\begin{equation}\label{5.53}
2^{5k/6} \| P_{k}(u^{5}) \|_{L_{t,x}^{6/5}([1, 2^{6(k - k_{0})}] \times \mathbf{R})} \lesssim \alpha(k)^{2} + 2^{-k/6} \alpha(k) \sum_{k_{0} \leq j \leq k} \alpha(j) 2^{(k - j)/10} 2^{j/6} (k - j + 1)^{5/8} 
\end{equation}

\begin{equation}\label{5.54}
 + 2^{-k/6} \alpha(k) \sum_{j \leq k_{0}} (k - k_{0})^{5/6} 2^{j/6}  \lesssim \alpha(k)^{2} + 2^{-k/4}.
\end{equation}

\noindent Also by $(\ref{5.44})$ and the proof of theorem $\ref{t5.3}$,

\begin{equation}\label{5.55}
 2^{5k/6} \| P_{> 2^{k}}(u^{5}) \|_{L_{t,x}^{6/5}([2^{-5k/2}, 1] \times \mathbf{R})} \lesssim \epsilon^{5}.
\end{equation}

\noindent Therefore, 

\begin{equation}\label{5.56}
(\ref{5.46}) \lesssim (\alpha(k)^{2} + 2^{-k/4}) \epsilon^{5}.
\end{equation}

\noindent Next,

\begin{equation}\label{5.57}
 \int_{0}^{2^{-5k/2}} \int_{1}^{2^{6(k - k_{0})}} \langle e^{(t - \tau) \partial_{x}^{3}} P_{k} \partial_{x}(u^{5}), P_{k} \partial_{x}(u^{5}) \rangle dt d\tau
\end{equation}

\begin{equation}\label{5.58}
 \lesssim 2^{2k} \| P_{k}(u^{5}) \|_{L_{t,x}^{1}([0, 2^{-5k/2}] \times \mathbf{R})} (\int_{2^{6(k - k_{0})}}^{\infty} \frac{1}{t^{1/3}} \| P_{k}(u^{5}) \|_{L_{x}^{1}} dt).
\end{equation}

\begin{equation}\label{5.59}
\aligned
 2^{k} \| P_{k}(u^{5}) \|_{L_{t,x}^{1}([T, 2T] \times \mathbf{R})} \\ \lesssim 2^{k} \| P_{> k - 5} u \|_{L_{x}^{\infty} L_{t}^{2}} \| u \|_{L_{t}^{\infty} L_{x}^{2}}^{1/2} \| P_{\leq 2^{k}} u \|_{L_{x}^{14/3} L_{t}^{14}}^{7/2} + 2^{k} \| P_{> 2^{k}} u \|_{L_{t,x}^{6}}^{9/2} \| u \|_{L_{t}^{\infty} L_{x}^{2}}^{1/2} \lesssim T^{1/4} 2^{k/4}.
\endaligned
\end{equation}

\noindent Therefore,

\begin{equation}\label{5.60}
 2^{k} \int_{2^{6(k - k_{0})}}^{\infty} \frac{1}{t^{1/3}} \| P_{2^{k}}(u^{5}) \|_{L_{x}^{1}} dt \lesssim 2^{-k/4}.
\end{equation}

\noindent By Holder's inequality,

\begin{equation}\label{5.61}
 2^{k} \| P_{k}(u^{5}) \|_{L_{t,x}^{1}([0, 2^{-5k/2}] \times \mathbf{R})} \lesssim 2^{k} \| u \|_{L_{t}^{\infty} L_{x}^{2}} \sum_{T < 2^{-5k/2}} T^{1/2} \| u \|_{L_{t,x}^{8}([T, 2T] \times \mathbf{R})}^{4} \lesssim 2^{-k/4}.
\end{equation}

\noindent Therefore,

\begin{equation}\label{5.62}
 \| P_{k} u(1) \|_{L^{2}}^{2} \lesssim \epsilon^{5} \alpha(k)^{2} + 2^{-k/2} + 2^{-k/8} \epsilon^{5} \alpha(k).
\end{equation}

\noindent Let $\beta(k)$ be another frequency envelope.

\begin{equation}\label{5.63}
\beta(k) = \sum_{j} 2^{-\frac{\delta}{2} |j - k|} \| P_{2^{k}} u(1) \|_{L^{2}(\mathbf{R})}.
\end{equation}

\begin{equation}\label{5.64}
\aligned
 \sum_{j} 2^{-\delta|j - k|/2} \| P_{j} u(1) \|_{L^{2}} \\ \lesssim \epsilon^{5/2} \sum_{j} 2^{-\delta|j - k|/2} \sum_{j_{1}} 2^{-\delta|j - j_{1}|} \| P_{j_{1}} u(1) \|_{L^{2}} + \sum_{j} 2^{-\delta |j - k|/2} 2^{-j/8}
\endaligned
\end{equation}

\noindent implies that for $\epsilon > 0$ sufficiently small,

\begin{equation}\label{5.65}
\beta(k) \lesssim \epsilon^{5/2} \beta(k) + 2^{-\delta k/16}.
\end{equation}

\noindent This implies after making the rescaling argument that

\begin{equation}\label{5.66}
 \mathcal M(2^{k}) \lesssim 2^{- \delta k/16}.
\end{equation}

\noindent Now suppose that for some $\sigma > 0$,

\begin{equation}\label{5.67}
 \mathcal M(A) \lesssim A^{-\sigma}.
\end{equation}

\begin{equation}\label{5.68}
 \| P_{> AT^{-1/3}} \partial_{x}(u^{5}) \|_{DU_{\partial_{x}^{3}}([T, 2T] \times \mathbf{R})} \lesssim \| P_{> A T^{-1/3}} u \|_{U_{\partial_{x}^{3}}([T, 2T] \times \mathbf{R})}^{5}
\end{equation}

\begin{equation}\label{5.69}
 + A^{-1/6} T^{1/18} \| \partial_{x} P_{> \frac{A}{32} T^{-1/3}} u \|_{L_{x}^{\infty} L_{t}^{2}([T, 2T] \times \mathbf{R})} \| P_{2^{k_{0}/2} T^{-1/3} < \cdot < \frac{A}{32} T^{-1/3}} u \|_{L_{x}^{24/5} L_{t}^{12}([T, 2T] \times \mathbf{R})}^{4} 
\end{equation}

\begin{equation}\label{5.70}
 + A^{-1/6} T^{1/18} \| \partial_{x} P_{> \frac{A}{32} T^{-1/3}} u \|_{L_{x}^{\infty} L_{t}^{2}([T, 2T] \times \mathbf{R})} \| P_{\leq 2^{k_{0}/2} T^{-1/3} } u \|_{L_{x}^{24/5} L_{t}^{12}([T, 2T] \times \mathbf{R})}^{4} \lesssim \mathcal S(\frac{A}{32}) \epsilon^{4} + 2^{\frac{k_{0}}{2}} A^{-1/6} \mathcal S(\frac{A}{32}).
\end{equation}

\noindent Since $\mathcal S(A) \lesssim \mathcal M(A) + \mathcal N(A)$,

\begin{equation}\label{5.71}
 \mathcal S(A) \lesssim A^{-\sigma} + \mathcal S(\frac{A}{32}) \epsilon^{4} + 2^{k_{0}/12} A^{-1/6} \mathcal S(\frac{A}{32}),
\end{equation}

\noindent so for $A \geq 2^{k_{0}}$, starting from $\mathcal S(2^{k_{0}}) \leq \epsilon$, by induction, taking $\epsilon(\sigma) > 0$ sufficiently small, here it suffices to consider $0 \leq \sigma \leq 2$,

\begin{equation}\label{5.72}
 \mathcal S(A) \lesssim A^{-\sigma}
\end{equation}

\noindent which in turn implies

\begin{equation}\label{5.73}
 \mathcal N(A) \lesssim A^{-\sigma}.
\end{equation}

\noindent Now we again use the no - waste Duhamel formula $(\ref{2.17.1})$.

\begin{equation}\label{5.74}
 \| P_{N} u(1) \|_{L^{2}} \lesssim \sum_{k \geq 0} \| P_{N} (u^{5}) \|_{L_{x}^{1} L_{t}^{2}([2^{k}, 2^{k + 1}] \times \mathbf{R})} \lesssim \sum_{k \geq 0} \| P_{> \frac{N}{32}} u \|_{L_{x}^{5} L_{t}^{10}([2^{k}, 2^{k + 1}] \times \mathbf{R})}^{5}
\end{equation}

\begin{equation}\label{5.75}
  + \sum_{k \geq 0} \| P_{> \frac{N}{32}} u \|_{L_{x}^{\infty} L_{t}^{2}([2^{k}, 2^{k + 1}] \times \mathbf{R})} \| P_{\leq N} u \|_{L_{x}^{4} L_{t}^{\infty}([2^{k}, 2^{k + 1}] \times \mathbf{R})}^{4}
\end{equation}

\begin{equation}\label{5.76}
 \lesssim (N^{-5 \sigma} + N^{-1 - \sigma}) \sum_{k \geq 0} 2^{-k \sigma} \lesssim N^{-5 \sigma} + N^{-1 - \sigma}.
\end{equation}

\noindent Iterating this argument finitely many times, this proves that $u(1) \in H^{1}$. This completes the proof of theorem $\ref{t5.1}$. $\Box$

\section{Rapid double cascade}
\begin{theorem}\label{t6.1}
 There does not exist a minimal mass blowup solution to the mass - critical gKdV in the form of a rapid double cascade.
\end{theorem}

\noindent \emph{Proof:} Let $t_{0} = t_{0}(T)$, where $t_{0}(T)$ is given by $(\ref{4.10})$. Let

\begin{equation}\label{6.1}
 u_{0}^{n} = \frac{1}{N(t_{0})^{1/2}} u(t_{0}, \frac{x + x(t_{0})}{N(t_{0})}).
\end{equation}

\noindent By concentration compactness $u_{0}^{n}$ has a subsequence that converges in $L^{2}$ to $u_{0} \in L^{2}$, and $u_{0}$ is the initial data for a minimal mass blowup solution to the mKdV on a maximal interval $I$, $N(0) = 1$, $N(t) \geq 1$ on $I$, and 

\begin{equation}\label{6.2}
 \int_{I} N(t)^{2} dt \lesssim C.
\end{equation}

\noindent Since $N(t) \geq 1$ this implies $|I| \lesssim C$, and also 

\begin{equation}\label{6.3}
\lim_{t \nearrow \sup(I)} N(t) = \lim_{t \searrow \inf(I)} N(t) = +\infty.
\end{equation}

\noindent Since $|x'(t)| \lesssim N(t)^{2}$ and $x(0) = 0$, $|x(t)| \lesssim C$ on $I$. Now define a Morawetz potential. Let $\psi \in C^{\infty}(\mathbf{R})$, $\psi$ is an odd function, $\psi(x) = x$ for $0 \leq x \leq 1$, $\psi(x) = \frac{3}{2}$ for $x > 1$. Also let

\begin{equation}\label{6.4}
 0 \leq \phi(x) = \psi'(x).
\end{equation}

\noindent For some $0 < R < \infty$ let

\begin{equation}\label{6.5}
 M(t) = R \int \psi(\frac{x}{R}) u(t,x)^{2} dx.
\end{equation}

\noindent For any $R > 0$, $N(t) \nearrow \infty$, as $t \rightarrow \sup(I), \inf(I)$, there exists $t_{+}$ sufficiently close to $\sup(I)$, $t_{-}$ sufficienly close to $\inf(I)$, such that

\begin{equation}\label{6.6}
 R \int \psi(\frac{x}{R}) u(t_{\pm},x)^{2} dx \lesssim C.
\end{equation}

\noindent Taking a derivative in time,

\begin{equation}\label{6.7}
 \frac{d}{dt} M(t) = -\int \phi(\frac{x}{R}) [3 u_{x}^{2} + \frac{5}{3} u^{6}] dx + O(\frac{1}{R^{2}}) \| u(t) \|_{L_{x}^{2}(\mathbf{R})}^{2}.
\end{equation}

\begin{equation}\label{6.8}
\frac{1}{R^{2}} \int_{I} \| u(t) \|_{L_{x}^{2}(\mathbf{R})}^{2} dt \lesssim \frac{C}{R^{2}}.
\end{equation}

\noindent Therefore, for any $R > 1$

\begin{equation}\label{6.9}
 \int_{I} \int_{|x| \leq R} [3 u_{x}^{2} + \frac{5}{3} u^{6}] dx dt \lesssim C.
\end{equation}

\noindent This bound is uniform in $R$, so in particular

\begin{equation}\label{6.10}
 \int_{I} \int [3 u_{x}^{2} + \frac{5}{3} u^{6}] dx dt \lesssim C.
\end{equation}

\noindent $(\ref{4.7})$ implies $|I| \gtrsim 1$. This in turn implies that there exists a $t \in I$ such that

\begin{equation}\label{6.11}
 \int [3 u_{x}^{2} + \frac{5}{3} u^{6}] dx \lesssim C.
\end{equation}

\noindent Conservation of energy then implies $E(u(t)) = E(u(0)) \lesssim C$ for all $t \in I$, which contradicts $N(t) \rightarrow +\infty$ as $t \rightarrow \sup(I)$ or $\inf(I)$. $\Box$

\section{Quasi - soliton}
\noindent Let 

\begin{equation}\label{7.1}
R(T) = C(\int_{0}^{T} N(t)^{2} dt)
\end{equation}

\noindent for some fixed constant $C$ such that $|x'(t)| \leq \frac{C}{2} N(t)^{2}$. $(\ref{4.15})$ implies

\begin{equation}\label{7.2}
 \sup_{t \in [0, T]} \int_{|x| \geq R(T)} u(t,x)^{2} dx \rightarrow 0,
\end{equation}

\noindent as $T \rightarrow \infty$. Once again let

\begin{equation}\label{7.3}
 M(t) = R \int \psi(\frac{x}{R}) u(t,x)^{2} dx.
\end{equation}

\begin{equation}\label{7.4}
 M(T) - M(0) \lesssim R.
\end{equation}

\begin{equation}\label{7.5}
\dot{M}(t) = -\int \phi(\frac{x}{R}) [3 u_{x}^{2} + \frac{5}{3} u^{6}] dx + \frac{1}{R^{2}} \int \phi''(\frac{x}{R}) u^{2} dx.
\end{equation}

\noindent For any $t_{0} \in [0, T]$,

\begin{equation}\label{7.6}
 \int_{t_{0}}^{t_{0} + \frac{\delta}{N(t_{0})^{3}}} 1 dt = \frac{\delta}{N(t_{0})^{3}} \lesssim (\int_{t_{0}}^{t_{0} + \frac{\delta}{N(t_{0})^{3}}} N(t)^{2} dt).
\end{equation}

\noindent This implies that since $N(0) = 1$,

\begin{equation}\label{7.7}
 \int_{0}^{T} 1 dt \lesssim (\int_{0}^{T} N(t)^{2} dt)^{3}
\end{equation}

\noindent which implies

\begin{equation}\label{7.8}
 \frac{1}{R^{2}} \int_{0}^{T} \int u(t,x)^{2} dx dt \lesssim \int_{0}^{T} N(t)^{2} dt.
\end{equation}

\noindent Fix $\mathcal J > 0$ large.

\begin{lemma}\label{l7.1}
There exists $I(T) \subset [0, T]$ with

\begin{equation}\label{7.9}
 \int_{I} N(t)^{3} dt = \mathcal J, \hspace{5mm} \int_{I} \int_{|x| \leq R(T)} [3 u_{x}^{2} + \frac{5}{3} u^{6}] dx dt \lesssim \int_{I} N(t)^{2} dt.
\end{equation}

\noindent The constant is uniform in $T$. 
\end{lemma}

\noindent Take $[0, T]$ such that

\begin{equation}\label{7.10}
\int_{0}^{T} N(t)^{3} dt = K \mathcal J
\end{equation}

\noindent for some integer $K$. Partition $[0, T]$ into intervals $I_{j}$. 

\begin{equation}\label{7.11}
 \sum_{j} \int_{I_{j}} \int_{|x| \leq R(T)} [3 u_{x}^{2} + \frac{5}{3} u^{6}] dx dt \lesssim \sum_{j} \int_{I_{j}} N(t)^{2} dt.
\end{equation}

\noindent Therefore there exists one $j$ such that

\begin{equation}\label{7.12}
 \int_{I_{j}} \int_{|x| \leq R(T)} [3 u_{x}^{2} + \frac{5}{3} u^{6}] dx dt \lesssim \int_{I_{j}} N(t)^{2} dt.
\end{equation}

\begin{lemma}\label{l7.2}
There exists $t_{0}(T) \in I(T)$ with

\begin{equation}\label{7.13}
 N(t_{0}) \lesssim (\frac{1}{\mathcal J} \int_{I} N(t)^{2} dt)^{-1},
\end{equation}

\begin{equation}\label{7.14}
\int_{|x| \leq R(T)} [3 u_{x}^{2} + \frac{5}{3} u^{6}] dx \lesssim N(t_{0})^{2}.
\end{equation}
\end{lemma}

\noindent \emph{Proof:} Suppose that for every $t$ with $N(t) \leq 10 (\frac{1}{\mathcal J} \int N(t)^{2} dt)^{-1}$,

\begin{equation}\label{7.15}
 \inf_{t \in J} \int_{|x| \leq R} [3 u_{x}^{2} + \frac{5}{3} u^{6}] dx >> N(t)^{2}.
\end{equation}

\noindent The contribution of these $N(t)$'s to $\int N(t)^{2} dt$ is small.

\begin{equation}\label{7.16}
 \int_{N(t) \geq 10 (\frac{1}{\mathcal J} \int N(t)^{2} dt)^{-1}} N(t)^{2} \leq \frac{1}{10(\frac{1}{\mathcal J} \int_{I} N(t)^{2} dt)^{-1}} \int_{I} N(t)^{3} dt \leq \frac{1}{10} (\int_{I} N(t)^{2} dt).
\end{equation}

\noindent Therefore $(\ref{7.15})$ implies

\begin{equation}\label{7.17}
\int_{I} \int_{|x| \leq R(T)} [3 u_{x}^{2} + \frac{5}{3} u^{6}] dx dt >> \int_{I(T)} N(t)^{2} dt,
\end{equation}

\noindent which contradicts $(\ref{7.12})$. $\Box$\vspace{5mm}

\noindent The sequence

\begin{equation}\label{7.18}
 \chi(\frac{x}{R(T)}) \frac{1}{N(t_{0}(T))^{1/2}} u(\frac{x - x(t_{0}(T))}{N(t_{0}(T))})
\end{equation}

\noindent has a subsequence that converges in $L^{2}$ to $u_{0} \in H^{1}$, $E(u_{0}) \lesssim 1$, and $u_{0}$ is the initial data for a minimal mass blowup solution to the mKdV problem.\vspace{5mm}

\noindent Moreover there exists an interval $I$, $0 \in I$, $\int_{I} N(t)^{3} dt = \mathcal J$ with

\begin{equation}\label{7.19}
 \int_{I} N(t)^{2} dt \lesssim \int_{I} N(t)^{3} dt \sim \mathcal J.
\end{equation}

\noindent By Holder's inequality,

\begin{equation}\label{7.20}
 \mathcal J^{2} \sim (\int_{I} N(t)^{3} dt)^{2} \lesssim (\int_{I} N(t)^{2} dt)(\int_{I} N(t)^{4} dt).
\end{equation}

\noindent This implies that

\begin{equation}\label{7.21}
 \int_{I} N(t)^{4} dt \gtrsim \mathcal J.
\end{equation}

\begin{theorem}[No quasi - soliton]\label{t7.3}
There does not exist a minimal mass blowup solution to $(\ref{1.1})$ satisfying $(\ref{7.19})$, $(\ref{7.21})$, $E(u(0)) \lesssim 1$ for $\mathcal J$ sufficiently large.
\end{theorem}

\noindent This theorem precludes the final minimal mass blowup solution since $\int N(t)^{3} dt$ is a scale invariant quantity and $(\ref{4.7})$ implies that $\int_{I} N(t)^{3} dt = +\infty$.\vspace{5mm}

\noindent \emph{Proof of theorem $\ref{t7.3}$:} We follow \cite{KwSh}, \cite{D3}, and especially \cite{TTao} to define an interaction Morawetz estimate. Recall $(\ref{3.1})$ - $(\ref{3.6})$. Define large constants $R$, $R_{1}$, $R_{1} << R$. Let $\chi_{a} \in C_{0}^{\infty}(\mathbf{R})$ be an even function, $\chi_{a} = 1$ for $|x| \leq a$, $\chi_{a} = 0$ for $|x| \geq a + R_{1}$, $a \geq R$. Let

\begin{equation}\label{7.22}
 \phi(x,y) = \frac{1}{R^{2}} \int_{R}^{2R} \int \chi_{a}(x - t) \chi_{a}(y - t) dt da.
\end{equation}

\begin{equation}\label{7.23}
 \phi(x,y) = \frac{1}{R^{2}} \int_{R}^{2R} \int \chi_{a}(x - y - t) \chi_{a}(t) dt da = \phi(x - y) = \frac{1}{R^{2}} \int_{R}^{2R} \int \chi_{a}(y + t - x) \chi_{a}(-t) dt = \phi(y - x).
\end{equation}

\noindent Then let

\begin{equation}\label{7.24}
 \psi(x - y) = \int_{0}^{x - y} \phi(t) dt.
\end{equation}

\noindent Now we produce an interaction Morawetz estimate. Let

\begin{equation}\label{7.25}
 M(t) = R \int \int \psi(\frac{(x - y) \tilde{N}(t)}{R}) \rho(t,y)^{2} e(t,x) dx dy.
\end{equation}

\noindent $\tilde{N}(t)$ is a quantity, $\tilde{N}(t) \leq N(t)$, that will be defined shortly.

\begin{equation}\label{7.26}
 \dot{M}(t) = \tilde{N}(t) \int \int \phi(\frac{(x - y) \tilde{N}(t)}{R}) [- \rho(t,y) k(t,x) + j(t,y) e(t,x)] dx dy 
\end{equation}

\begin{equation}\label{7.27}
 + \frac{\tilde{N}(t)^{3}}{R^{2}} \int \int \rho(t,y)^{2} e(t,x)^{2} dx dy 
\end{equation}

\begin{equation}\label{7.28}
+ \int \int \frac{\tilde{N}'(t) (x - y)}{R} \phi(\frac{(x - y)\tilde{N}(t)}{R}) u(t,y)^{2} [\frac{1}{2} u_{x}^{2} + u^{6}] dx dy.
\end{equation}

\begin{equation}\label{7.29}
 (\ref{7.26}) = - \tilde{N}(t) \int \int \int \chi(x \tilde{N}(t) - s) \chi(y \tilde{N}(t) - s) [u(t,y)^{2} (\frac{3}{2} u_{xx}^{2} + 2 u_{x}^{2} u^{4} + \frac{1}{2} u^{10})] dx dy
\end{equation}

\begin{equation}\label{7.30}
 + \tilde{N}(t) \int \int \int \chi(x \tilde{N}(t) - s) \chi(y \tilde{N}(t) - s) [3 u_{y}^{2} + \frac{5}{3} u^{6}] [\frac{1}{2} u_{x}^{2} + \frac{1}{6} u^{6}] dx dy
\end{equation}

\noindent Let $\tilde{\chi} = 1$ on $[a, a + R_{1}]$ and $0$ elsewhere. We will suppress the $a$ for the moment and take $\chi_{a} = \chi$ for some $a$.

\begin{equation}\label{7.31}
 \int \chi^{2} u_{xx}^{2} dx = \int \chi u_{xx} [\partial_{xx}(\chi u) - 2 \chi_{x} u_{x} - \chi_{xx} u] dx
\end{equation}

\begin{equation}\label{7.32}
 = \int \chi u_{xx} \partial_{xx}(\chi u) dx - \int \chi_{x} \chi \partial_{x}(u_{x}^{2}) dx - \int \chi_{xx} \chi u_{xx} u dx
\end{equation}

\begin{equation}\label{7.33}
\aligned
 = \int \partial_{xx}(\chi u)^{2} dx - 2 \int \chi_{x} u_{x} \partial_{xx}(\chi u) dx - \int \chi_{xx} u \cdot \partial_{xx}(\chi u) dx \\ + \int \frac{1}{2} \partial_{xx}(\chi^{2}) u^{2} dx + \int \chi_{xx} \chi u_{x}^{2} dx - \frac{1}{2} \int \partial_{xx}(\chi_{xx} \chi) u^{2} dx
\endaligned
\end{equation}

\begin{equation}\label{7.34}
 = \int \partial_{xx}(\chi u)^{2} dx + \frac{1}{R_{1}^{2}} \int O(u_{x}^{2} \tilde{\chi}^{2}) dx + \frac{1}{R_{1}^{4}} \int O(\tilde{\chi}^{2} u^{2}) dx.
\end{equation}

\noindent Next,

\begin{equation}\label{7.35}
 \int \chi^{2} u_{x}^{2} = \int \chi u_{x} \partial_{x}(\chi u) - \int \chi u_{x} \chi_{x} u
\end{equation}

\begin{equation}\label{7.36}
 = \int \partial_{x}(\chi u)^{2} + \frac{1}{4} \int \partial_{xx}(\chi^{2}) u^{2} - \int \chi_{x} u \partial_{x}(\chi u) = \int \partial_{x}(\chi u)^{2} + \frac{1}{R_{1}^{2}} \int \tilde{\chi}^{2} u^{2}.
\end{equation}

\noindent Next,

\begin{equation}\label{7.37}
 \int \chi^{2} u_{x}^{2} u^{4} = \int \partial_{x}(\chi u) \chi u_{x} u^{4} - \int \chi_{x} u \chi u_{x} u^{4}
\end{equation}

\begin{equation}\label{7.38}
 = \int \partial_{x}(\chi u)^{2} u^{4} - \int \chi_{x} u^{5} \partial_{x}(\chi u) + \frac{1}{2} \int \partial_{xx}(\chi^{2}) u^{6}
\end{equation}

\begin{equation}\label{7.39}
 \geq \int \partial_{x}(\chi u)^{2} (\chi u)^{4} + \frac{1}{R_{1}^{2}} \int \tilde{\chi}^{2} u^{6}.
\end{equation}

\noindent Finally,

\begin{equation}\label{7.40}
 \int \chi^{2} u^{6} dx = \int (\chi u)^{6} dx + \int (1 - \chi^{4}) (\chi u)^{2} u^{4}.
\end{equation}

\noindent From \cite{TTao}, if $v = \chi_{a} u$,

\begin{equation}\label{7.41}
 \frac{3}{2}(\int v^{2})(\int v_{xx}^{2}) - \frac{3}{2} (\int v_{x}^{2})^{2} + 2 (\int v_{x}^{2} v^{4})(\int v^{2}) + \frac{1}{2} (\int v^{10})(\int v^{2})
\end{equation}

\begin{equation}\label{7.42}
 - \frac{4}{3} (\int v^{6})(\int v_{x}^{2}) - \frac{1}{2} (\int v^{6} )^{2} > 0.
\end{equation}

\noindent Next, for $R$ sufficiently large, by Holder's inequality,

\begin{equation}\label{7.43}
 \frac{2}{9R} \int (\int \chi_{a}^{6}(\frac{x \tilde{N}(t)}{R} - s) u(t,x)^{6} dx)^{2} ds \gtrsim (\int_{|x - x(t)| \leq \frac{C_{0}}{N(t)}} u(t,x)^{6} dx) \gtrsim N(t)^{4},
\end{equation}

\noindent uniformly in $a$. Now we estimate the contribution of the errors. Let $1_{A}(x)$ be the indicator function of a set $A$.

\begin{equation}\label{7.44}
 \frac{1}{R} \int_{R}^{2R} 1_{[a, a + R_{1}]} da \leq \frac{R_{1}}{R} 1_{[R, 3R]}.
\end{equation}

\noindent By Holders inequality, Sobolev embedding, and $(\ref{4.7})$,

\begin{equation}\label{7.45}
\aligned
 \| u \|_{L_{t,x}^{6}([t_{0}, t_{0} + \frac{\delta}{N(t_{0})^{3}}] \times \mathbf{R})}^{6} \lesssim \| u_{\geq N(t_{0})} \|_{L_{t,x}^{6}([t_{0}, t_{0} + \frac{\delta}{N(t_{0})^{3}}] \times \mathbf{R})}^{6} + \| u_{\leq N(t_{0})} \|_{L_{t,x}^{6}([t_{0}, t_{0} + \frac{\delta}{N(t_{0})^{3}}] \times \mathbf{R})}^{6} \\ \lesssim \frac{1}{N(t_{0})} + N(t_{0})^{2} \frac{1}{N(t_{0})^{3}} \sim \int_{t_{0}}^{t_{0} + \frac{\delta}{N(t_{0})^{3}}} N(t)^{2} dt.
 \endaligned
\end{equation}

\noindent Therefore, by conservation of energy

\begin{equation}\label{7.46}
 \frac{R_{1}}{R} \int_{I} \int \int u(t,y)^{6} u_{x}^{2} dx dy dt \lesssim \frac{R_{1}}{R} \int_{I} N(t)^{2} dt.
\end{equation}

\noindent Next, let $J_{l} = [t_{0}, t_{0} + \frac{\delta}{N(t_{0})^{3}}]$.

\begin{equation}\label{7.47}
 \frac{R_{1}}{R} \frac{\tilde{N}(t)^{3}}{R_{1}^{2}} \int_{J_{l}} (\int_{|x - y| \lesssim \frac{R}{\tilde{N}(t)}} u_{x}^{2} u^{2}) 
\end{equation}

\begin{equation}\label{7.48}
\lesssim  \frac{\tilde{N}(t)^{3}}{R_{1} R} \| u_{x}(u^{2}) \|_{L_{t,x}^{2}(J_{l} \times \mathbf{R})} \| u \|_{L_{t}^{\infty} L_{x}^{2}(J_{l} \times \mathbf{R})} \frac{R}{\tilde{N(t_{0})}} \frac{1}{N(t_{0})^{3/2}} \lesssim \frac{1}{R_{1}} \frac{\tilde{N}(t_{0})^{2}}{N(t_{0})^{3/2}}.
\end{equation}

\noindent The last inequality follows from conservation of energy, Holder's inequality, and 

\begin{equation}\label{7.49}
 \| u \|_{L_{x}^{4} L_{t}^{\infty}(J_{l} \times \mathbf{R})} \lesssim \| \partial_{x} u \|_{S^{0}(J_{l} \times \mathbf{R})}^{1/4} \| u \|_{S^{0}(I \times \mathbf{R})}^{3/4},
\end{equation}

\begin{equation}\label{7.50}
 \| u \|_{L_{x}^{\infty} L_{t}^{2}(J_{l} \times \mathbf{R})} \lesssim \| u \|_{S^{0}(I \times \mathbf{R})}.
\end{equation}

\noindent Next, by conservation of mass

\begin{equation}\label{7.51}
 \frac{\tilde{N}(t_{0})^{4}}{R R_{1}^{3}} \int_{J_{l}} \int_{|x - y| \sim \frac{R}{\tilde{N}(J)}} u(t,x)^{2} u(t,y)^{2} dx dy \lesssim \frac{\tilde{N}(t_{0})}{R R_{1}^{3}}.
\end{equation}

\noindent Finally, by conservation of mass and $(\ref{7.45})$

\begin{equation}\label{7.52}
 \int_{I} \frac{\tilde{N}(t)^{3}}{R_{1} R} \int u(t,x)^{2} u(t,y)^{6} dx dy dt \lesssim \frac{1}{R_{1} R} \int \tilde{N}(t)^{2} N(t)^{2} dt.
\end{equation}

\noindent This takes care of the error terms in $(\ref{7.34})$, $(\ref{7.36})$, $(\ref{7.39})$, and $(\ref{7.40})$. By the fundamental theorem of calculus and the above computations, taking say $R_{1} = R^{1/2}$,

\begin{equation}\label{7.53}
 \int_{I} N(t)^{4} \tilde{N}(t) dt \lesssim \eta(R) \int_{I} N(t)^{2} \tilde{N}(t) dt + R \int \frac{|\tilde{N}'(t)|}{\tilde{N}(t)} \int_{|x - y| \lesssim \frac{R}{\tilde{N}(J)}} u_{x}^{2} u^{2} dx dy dt,
\end{equation}

\noindent where $\eta(R) \rightarrow 0$ as $R \rightarrow \infty$. Now choose $\tilde{N}(t) = N(t)$ for $N(t) \leq \alpha$ and $\tilde{N}(t) = \alpha$ for $N(t) \geq \alpha$, $\alpha > 0$ some small fixed constant.

\begin{equation}\label{7.54}
 R \int_{J_{l}} \tilde{N}(t_{0})^{3} \int_{|x - y| \lesssim \frac{R}{\tilde{N}(t_{0})}} u_{x}^{2} u^{2} dx dy dt \lesssim \frac{R^{3/2}}{N(t_{0})^{3/2}} \tilde{N}(t_{0})^{2}  \sim R^{3/2} \int_{J_{l}} \tilde{N}(t)^{2} N(t)^{3/2} dt.
\end{equation}

\noindent Since $\tilde{N}(t) \leq \alpha$, $\tilde{N}(t) \leq N(t)$,

\begin{equation}\label{7.55}
 R^{3/2} \int_{I} \tilde{N}(t)^{2} N(t)^{3/2} \lesssim R^{3/2} \alpha^{3/2} \mathcal J.
\end{equation}

\begin{equation}\label{7.56}
 \eta(R) \int_{I} \tilde{N}(t) N(t)^{2} dt \lesssim \eta(R) \alpha \mathcal J.
\end{equation}

\noindent Next,

\begin{equation}\label{7.57}
 \int_{t : N(t) \leq \alpha} N(t)^{4} dt \leq \alpha^{2} \mathcal J.
\end{equation}

\noindent Since $\int_{I} N(t)^{4} dt \gtrsim \mathcal J$, by the fundamental theorem of calculus and the error estimates,

\begin{equation}\label{7.58}
 \alpha \mathcal J \lesssim \alpha \int_{I} N(t)^{4} dt \lesssim R + \eta(R) \alpha \mathcal J + R^{3/2} \alpha^{3/2} \mathcal J.
\end{equation}

\noindent Choose $\alpha(R)$ sufficiently small so that $\alpha^{3/2} R^{3/2} << \eta(R)$. Then for $\mathcal J$ sufficiently large, we have a contradiction. $\Box$\vspace{5mm}

\nocite*
\bibliographystyle{plain}

\end{document}